\documentclass[a4paper,12pt]{article}
\usepackage{geometry,latexsym,amssymb,amsmath,amsthm,color}
\usepackage[latin5]{inputenc}
\usepackage{enumerate}
\usepackage{enumitem}
\usepackage[T1]{fontenc}
\usepackage{authblk}
\usepackage[all]{xy}
\usepackage{palatino}
\usepackage{indentfirst}
\usepackage{setspace}
\newtheorem{example}{Example}[section]
\newtheorem{Def}[example]{Definition}
\newtheorem{Exam}[example]{Example}
\newtheorem{Prop}[example]{Proposition}
\newtheorem{Theo}[example]{Theorem}

\newtheorem{Cor}[example]{Corollary}
\newenvironment{Prf}{{\bf Proof:} }{\hfill $\Box$
	\mbox{}}
\def\wtilde{\widetilde}

\begin{document}
\title{Covering groupoids of categorical groups}

\author[a]{Osman Mucuk\thanks{E-mail : mucuk@erciyes.edu.tr}}
\author[b]{Tunçar Şahan\thanks{E-mail : tuncarsahan@aksaray.edu.tr}}
\affil[a]{\small{Department of Mathematics, Erciyes University, Kayseri, TURKEY}}
\affil[b]{\small{Department of Mathematics, Aksaray University, Aksaray, TURKEY}}

\date{}

\maketitle
\begin{abstract}
If $X$ is a topological group, then its fundamental groupoid
$\pi_1X$ is a  group-groupoid which is  a group object in the
category of groupoids. Further if $X$ is a path connected
topological group which has a simply connected cover, then the
category of covering spaces   of $X$ and the category of covering
groupoids of $\pi_1X$ are equivalent.  In this paper we prove  that
if $(X,x_0)$ is an $H$-group, then the fundamental groupoid $\pi_1X$
is a categorical group.  This enable us to prove that the category
of the covering spaces of an $H$-group $(X,x_0)$ is equivalent to
the category of covering groupoid of the categorical group $\pi_1X$.
\end{abstract}

\\ {\bf Classification:} 22A05, 55M99, 55R15
	
\section{Introduction}

Covering spaces are studied in algebraic topology, but they have
important applications
in many other branches of mathematics including differential topology,
the theory of topological  groups and the theory of Riemann surfaces.

One of the ways of expressing the algebraic content of the theory
of covering spaces is using groupoids and the fundamental
groupoids.  The latter functor gives an equivalence of categories
between the category of covering spaces of a reasonably nice space
$X$ and the category of groupoid covering morphisms of $\pi_1 X$.

If $X$ is a connected topological group with identity $e$, $
p\colon (\widetilde{X},\widetilde{e})\rightarrow (X,e)$ is a
covering map of pointed spaces such that $\widetilde{X}$ is simply
connected, then $\wtilde X$ becomes a topological group with
identity $\widetilde{e}$ such that  $p$ is a morphism of
topological groups \cite{Ch,Rot}.

The problem of universal covers of
non-connected topological groups was first studied by Taylor in
~\cite{Ta}. He proved that a topological group $X$ determines an
obstruction class $k_X$ in $H^3(\pi_0(X),\pi_1(X,e))$, and that
the vanishing of $k_{X}$ is a necessary and sufficient condition
for the lifting of the group structure to a universal cover. In
~\cite{Mu1} an analogous algebraic result is given in terms of
crossed modules and group objects in the category of groupoids
(see also ~\cite{BM1} for a revised version, which generalizes
these results and shows the relation with the theory of
obstructions to extension  for groups).

For a topological group  $X$,  the fundamental groupoid
$\pi_1 X$  becomes  a group object in the category of groupoids
\cite{BS}.  This is also called an internal category in the
category of groups \cite{Ti}. This functor  gives an equality of
the category of the covering spaces of a topological group $X$
whose underlying space is locally nice  and the category  of the
covering groupoids of $\pi_1X$ \cite{BM1,Mu1}.

In this paper we prove  that if $(X,x_0)$ is an $H$-group, then
the fundamental groupoid $\pi_1X$ is a categorical group.  This enable us to prove that the category of the covering spaces of an $H$-group $(X,x_0)$ is equivalent to the category of covering groupoid of the categorical group $\pi_1X$.

\section{ Covering Spaces and $H$-groups} \label{Section1} We assume the usual
theory of covering maps. All spaces $X$ are assumed to be locally path
connected and semi locally 1-connected, so that each path
component of $X$  admits a simply connected cover.

Recall that a  covering map $p \colon \wtilde{X} \rightarrow X$ of
connected spaces is called  {\it universal}  if it covers every
cover of $X$ in the sense that if $q \colon \wtilde{Y} \rightarrow
X$ is another cover of $X$ then there   exists a map  $r \colon
\wtilde{X} \rightarrow \wtilde{Y}$ such that $p=qr$ (hence $r$
becomes a cover). A covering map $p \colon \wtilde{X} \rightarrow
X$ is called {\it simply connected} if $\wtilde{X}$ is simply
connected. Note that a simply connected cover is a universal
cover.
\begin{Def}
	{\em We call a  subset $U$ of $X$  {\em liftable} if it is open,
		path connected and $U$ lifts to each cover of $X$, that is, if $p
		\colon \wtilde{X} \rightarrow X$ is a covering map,  $\imath
		\colon U\rightarrow X$ is the inclusion map, and $\wtilde{x} \in
		\wtilde{X}$ satisfies $p(\wtilde x)=x \in U$, then there exists a
		map (necessarily unique) $\widehat{\imath} \colon U\rightarrow
		\wtilde{X}$ such that $p\widehat{\imath}=\imath$ and
		$\widehat{\imath}(x)=\wtilde x$.}
\end{Def}

It is easy to see that $U$ is liftable if and only if it is open,
path connected and for all $x\in U$, $\imath_{\star}\pi_{1}(U,x)$
is singleton, where $\pi_{1}(U,x)$ is the fundamental group of $U$
at the base point $x$ and  $\imath_{\star}$ is the morphism
$\pi_1(U,x)\rightarrow \pi_1(X,x)$ induced by the inclusion map
$\imath\colon U\rightarrow X$.  Remark that if $X$ is a semi
locally simply connected topological space, then each point $x\in
X$ has a liftable neighbourhood.

\begin{Def} {\em (~\cite{Rot}) \label{LiftLemma}
		Let $X$ be a topological space. Two covering maps $p\colon
		(\widetilde{X},\widetilde{x}_0)\rightarrow (X,x_0)$ and $q\colon
		(\widetilde{Y},\widetilde{y}_0)\rightarrow (X,x_0)$ are called
		equivalent if there is a homeomorphism $f\colon
		(\widetilde{X},\widetilde{x}_0)\rightarrow
		(\widetilde{Y},\widetilde{y}_0)$ such that $qf=p$.}
\end{Def}
\begin{Def} {\em Let $p\colon (\widetilde{X},\widetilde{x}_0)\rightarrow
		(X,x_0)$ be a covering map. The subgroup
		$p_{\star}(\pi_1(\widetilde{X},\widetilde{x}_0))$ of
		$\pi_1(X,x_0)$ is called {\em characteristic group} of $p$, where
		$p_{\star}$ is the morphism induced by $p$.}\end{Def}

The following result is known as Covering Homotopy Theorem
(Theorem 10.6 in ~\cite{Rot}).

\begin{Theo}\label{homotopy covering theorem}{\em  \label{CovHomTheo} Let
		$p\colon \widetilde{X} \rightarrow X$ be a covering
		map and $Z$ a connected space. Consider the commutative  diagram
		of continuous maps
		\begin{equation*}
			\xymatrix{
				Z \ar[d]_j \ar[r]^{\widetilde{f}}
				& \widetilde{X}  \ar[d]^p \\
				Z\times I \ar@{-->}[ur]^{\widetilde{F}} \ar[r]_>>>>>>F & X }
		\end{equation*}
		
		where $j\colon Z\rightarrow Z\times I, j(z)=(z,0)$ for all $z\in
		Z$. Then there is a unique continuous map $\widetilde{F}\colon
		Z\times I\rightarrow\widetilde{ X}$ such that $p\widetilde{F}=F$
		and $\widetilde{F} j=\widetilde{f}$ }\end{Theo}

As a corollary of this theorem if the maps $f,g\colon Z\rightarrow
X$ are homotopic, then their respective liftings $\widetilde{f}$
and  $\widetilde{g}$ are homotopic.  If $f\simeq g$, there is a
continuous map $F\colon Z\times I\rightarrow X$ such that
$F(z,0)=f(z)$ and $F(z,1)=g(z)$. So there is a continuous map
$\widetilde{F}\colon Z\times I\rightarrow \widetilde{X}$ as in
Theorem \ref{homotopy covering theorem}. Here
$p\widetilde{F}(z,0)=F(z,0)=f(z)$ and
$p\widetilde{F}(z,1)=F(z,1)=g(z)$. By the uniqueness of the
liftings  we have that $\widetilde{F}(z,0)=\widetilde{f}(z)$ and
$\widetilde{F}(z,1)=\widetilde{g}(z)$. Therefore
$\widetilde{f}\simeq \widetilde{g}$.

\begin{Def}\label{DefHgroup}{\em (\cite{Rot}, p.324)} {\em   A pointed space
		$(X,x_0)$ is called  an $H-group$ if there are continuous
		pointed maps \[ m \colon (X\times X,(x_0,x_0))\rightarrow (X,x_0),(x,x')\mapsto
		x\circ x'\]
		\[n\colon (X,x_0)\rightarrow (X,x_0)\]
		and pointed
		homotopies
		
		(i)  $m(1_X\times m)\simeq m(m\times 1_X)$
		
		(ii) $m~i_1\simeq 1_X\simeq m~i_2$
		
		(iii) $m(1_X,n)\simeq c\simeq m(n,1_X)$ \\ where $i_1,i_2\colon
		X\rightarrow X\times X $ are injections defined by
		$i_1(x)=(x,x_0)$ and $i_2(x)=(x_0,x)$, and  $c\colon X\rightarrow
		X$ is the constant map at $x_0$.}\end{Def}

In  an $H$-group $X$,  we denote $m(x,x')$ by $x\circ x'$ for
$x,x'\in X$ .

Let $(X,x_0)$ and $(Y,y_0)$ be $H$-groups. A continuous map
$f\colon (X,x_0)\rightarrow (Y,y_0)$ such that $f(x\circ
x')=f(x)\circ f(x')$ for $x,x'\in X$ is called a {\em morphism} of
$H$-groups. So we have a category of $H$-groups denoted by
$HGroups$.

\begin{Exam} {\em A topological group $X$ with identity $e$ is an
		$H$-group.  For  a topological group $X$,  the group operation
		\[ m \colon (X\times X,(e,e))\rightarrow (X,e)\] and the inverse
		map $n\colon (X,x_0)\rightarrow (X,x_0)$ are continuous; and
		
		(i) $m(1_X\times m)= m(m\times 1_X)$
		
		(ii) $m~ i_1=1_X= m~i_2$
		
		(iii) $m(1_X,n)= c= m(n,1_X)$,  where $c\colon X\rightarrow X$ is
		the constant map at $e$.}\end{Exam}

\begin{Theo}\label{TheoRotman}{\em (Theorem 11.9 in   \cite{Rot})}  If  $(X,x_0)$
	is a pointed space,
	then the loop space $\Omega(X,x_0)$ is an $H$-group.\end{Theo}

\begin{Def}{\em
		Let $(X,x_0)$ and $(Y,y_0)$ be $H$-groups and $U$ an open
		neighbourhood of  $x_0$ in $X$. A continuous map
		$\phi \colon (U,x_0)\rightarrow (Y,y_0)$ is called a {\em local  morphism}
		of $H$-groups if  $\phi(x\circ y)=\phi(x)\circ \phi(y)$ for
		$x,y\in U$ such that $x\circ y\in U$.}\end{Def}

\begin{Theo}\label{PropLift} Let $(X,x_0)$ and
	$(\wtilde{X},\widetilde{x}_0)$ be  $H$-groups  and $q \colon
	(\wtilde{X},\widetilde{x}_0)\rightarrow (X,x_0)$ a morphism of
	$H$-groups which is a covering map on the underlying spaces. Let
	$U$ be an open, path connected neighbourhood of $x_0$ in $X$ such
	that $U^{2}$ is contained in a liftable neighbourhood $V$ of $x_0$
	in
	$X$. Then the inclusion map $\imath\colon (U,x_0)\rightarrow (X,x_0)$
	lifts to a local morphism $\widehat{\imath} \colon
	(U,x_0)\rightarrow (\wtilde{X},\widetilde{x}_0)$ of H-groups.
\end{Theo}
\begin{Prf}
	Since $V$ lifts to $\wtilde{X}$, then $U$ lifts to $\wtilde{X}$ by
	$\widehat{\imath}\colon (U,x_0)\rightarrow
	(\wtilde{X},\widetilde{x}_0)$. We now prove that
	$\widehat{\imath}$ is a local morphism of $H$-groups. By the
	lifting theorem the map $\widehat{\imath}\colon (U,x_0)\rightarrow
	(\wtilde{X},\widetilde{x}_0)$ is continuous. We have to prove that
	$\widehat{\imath}\colon (U,x_0)\rightarrow
	(\wtilde{X},\widetilde{x}_0)$ preserves the multiplication. Let
	$x,y\in U$ such that $x\circ y\in U$. Let $a$ and $b$ be  paths
	from $x$ and $y$ to $x_0$ in $U$ respectively.  By the continuity
	of
	\[ m \colon (X\times X,(x_0,x_0))\rightarrow (X,x_0)\]
	$c=a\circ b$ defined by $c(t)=a(t)\circ b(t)$ for
	$t\in [0,1]$ is a  path from $x\circ y$ to $x_0$. Since
	$U^{2}\subseteq V$, the path $c$ is in $V$. So the paths $a,b$ and
	$c$ lift to $\wtilde{X}$. Let $\wtilde{a}$, $\wtilde{b}$ and
	$\wtilde{c}$ be  the liftings  of the paths $a$, $b$ and $c$  with  end points
	$x_0$ as chosen above.
	respectively. Since $q$ is a morphism of $H$-spaces, we have
	\[       q(\wtilde{c})=c=a\circ b= q(\wtilde{a})\circ q(\wtilde{b}).\]
	and
	\[     q(\wtilde{a}\circ \wtilde{b})=q(\wtilde{a})\circ q(\wtilde{b}).\]
	Since $\wtilde c$ and $\wtilde a\circ\wtilde b$ end at
	$\widetilde{x}_0\in \wtilde X$, by the  unique path lifting, we
	have that
	\[ \wtilde c=\wtilde{a}\circ\wtilde{b} \]
	By evaluating these paths at $0\in I$ we
	have that
	\[       \widehat{\imath}(x\circ y) =\widehat{\imath}(x)\circ
	\widehat{\imath}(y). \]
	Hence $\widehat{\imath}\colon (U,x_0)\rightarrow
	(\wtilde{X},\widetilde{x}_0)$ is a local morphism of $H$-groups.
\end{Prf}

\section{Covering Groupoids}
A {\em groupoid} $G$ on $O_G$ is a small category in which each
morphism is
an isomorphism. Thus $G$ has a set of morphisms,  a set $O_G$ of {\em objects}
together with functions
$s, t\colon G\rightarrow O_G$, $\epsilon \colon O_G \rightarrow G$
such that $s\epsilon=t\epsilon=1_{O_G}$, the identity map. The functions
$s$, $t$ are called {\em initial} and {\em final} point maps respectively.
If $a,b\in G$ and $t(a)=s(b)$, then the {\em product} or {\em composite}
$ba$ exists such that $s(ba)=s(a)$ and $t(ba)=t(b)$. Further,
this composite is associative, for $x\in O_G$ the element $\epsilon (x)$ denoted
by $1_x$ acts as the identity, and each element $a$ has an inverse $a^{-1}$ such
that $s(a^{-1})=t(a)$, $t(a^{-1})=s(a)$, $aa^{-1}=(\epsilon t)(a)$,
$a^{-1}a=(\epsilon s)(a)$. The map $G\rightarrow G$, $a\mapsto a^{-1}$,
is called the {\em inversion}.

So a group is a groupoid with only one object.

In a groupoid $G$
for $x,y\in O_G$ we write $G(x,y)$ for the set of all
morphisms with initial point $x$ and final point $y$. We say $G$ is
{\em connected} if for all $x,y \in O_G$, $G(x,y)$ is not empty and {\em simply
	connected} if $G(x,y)$ has only
one morphism. For
$x \in O_G$ we denote the star $\{a\in G\mid s(a)=x\} $ of $x$ by
$G_x$. The {\em object group} at $x$ is $G(x)=G(x,x)$.

Let $G$ and $H$ be groupoids. A {\em morphism} from $H$ to $G$ is
a pair of maps
$f\colon H\rightarrow G$ and $O_f\colon O_H\rightarrow O_G$
such that $s\circ f=O_f\circ s$, $t\circ f=O_f\circ t$ and
$f(ba)=f(b)f(a)$ for all $(a,b)\in H_{t}\times_{s} H$. For such a
morphism we simply write $f\colon H\rightarrow G$.

\begin{Def}{\em  Let $p\colon\wtilde G\rightarrow G$ be a morphism of groupoids.
		Then $p$ is
		called a {\em covering morphism} and $\widetilde{G}$ {\em a covering groupoid}
		of $G$ if for
		each $\wtilde x\in O_{\wtilde G}$ the restriction of $p$ \[p_x\colon
		(\wtilde{G})_{\wtilde x}\rightarrow G_{p(\wtilde x)}\]  is
		bijective. A covering morphism $p\colon \widetilde{G}\rightarrow
		G$ is called {\em connected } if both $\widetilde{G}$ and $G$ are
		connected.}\end{Def}

A connected covering morphism $p\colon\wtilde G\rightarrow G$
is called {\em universal} if $\wtilde G$ covers every cover of $G$, i.e. if
for every covering morphism $q\colon \widetilde{H}\rightarrow G$ there is a
unique
morphism of groupoids $\widetilde{p}\colon \wtilde G\rightarrow \widetilde{H}$
such that $q\widetilde{p}=p$
(and hence $\widetilde{p}$ is also a covering morphism), this is equivalent to
that for
$\wtilde{x}, \wtilde{y}\in O_{\wtilde G}$ the set $\wtilde{G}(\wtilde x, \wtilde
y)$
has not more than one element.

A group homomorphism $f\colon G\rightarrow H$ is a covering
morphism  if and only if it is an isomorphism.

For any groupoid morphism $p\colon \widetilde{G}\rightarrow G$ and
an object $\widetilde{x}$ of $\widetilde{G}$ we call  the subgroup
$p(\widetilde{G}(\widetilde{x}))$ of $G(p\widetilde{x})$ the {\em
	characteristic group} of $p$ at $\widetilde{x}$.

\begin{Exam}{\em \cite{Br1}} \label{fundgpd}{\em  If $p\colon
		\widetilde{X}\rightarrow X$ is a covering map of topological
		spaces, then the induced fundamental groupoid morphism
		$\pi_1(p)\colon \pi_1(\widetilde{X})\rightarrow \pi_1(X)$ is a
		covering morphism of groupoids. }\end{Exam}

\begin{Def} Let $p\colon\tilde G\rightarrow G$ be a covering morphism of
	groupoids and $q\colon H\rightarrow G$ a morphism of groupoids. If
	there exists a unique morphism $\tilde q\colon H\rightarrow \tilde
	G$ such that $q=p \tilde q$ we just say q lifts to $\tilde q$ by
	p.
\end{Def}
We call the following theorem from Brown [1] which is an important
result to have the lifting maps on covering groupoids.
\begin{Theo} \label{Theolifted}
	Let $p\colon\wtilde G\rightarrow G$ be a covering morphism of
	groupoids, $x\in O_G$ and $\wtilde x\in O_{\wtilde G}$ such that
	$p(\wtilde x)=x$. Let $q\colon K\rightarrow G$ be a morphism of
	groupoids such that $K$ is connected and $z\in O_K$ such that
	$q(z)=x$. Then the morphism $q\colon K\rightarrow G$ uniquely
	lifts to a morphism $\wtilde q\colon K\rightarrow \wtilde G$ such
	that $\wtilde q(z)=\wtilde x$ if and only if $q[K(z)]\subseteq
	p[\wtilde G(\wtilde x)]$, where $(z)$ and $\wtilde G(\wtilde x)$
	are the object groups.
\end{Theo}

\begin{Cor} \label{CorLiftcov} Let $p\colon
	(\widetilde{G}(\widetilde{x})\rightarrow (G,x)$ and
	$q\colon (\widetilde{H},\widetilde{z})\rightarrow (G,x)$ be connected
	covering morphisms with characteristic groups $C$ and $D$
	respectively. If $C\subseteq D$, then there is a unique covering
	morphism $r\colon (\widetilde{G},\widetilde{x})\longrightarrow
	(\widetilde{H},\widetilde{z})$  such that $p=qr$. If $C=D$, then $r$
	is an isomorphism.\end{Cor}

\section{Homotopies of functors and categorical groups}

In this section we prove that the functors are homotopic if and only if they are naturally isomorph.  For the homotopies of functors we first need  the following fact whose proofs are straightforward [1].
\begin{Prop}{\em Let $\mathcal{C}$, $\mathcal{D}$ and
		$\mathcal{E}$ be categories and  $F\colon \mathcal{C}\times
		\mathcal{D}\rightarrow \mathcal{E}$  a functor.  Then for $x\in
		Ob(C)$ and $y\in Ob(D)$ we have the induced functors
		
		\[F(x,-)\colon \mathcal{D}\rightarrow \mathcal{E}\]
		\[F(-,y)\colon \mathcal{C}\rightarrow \mathcal{E}.\]
	}\end{Prop}
	
	We  write  $\mathcal{J}$ for the simply connected groupoid whose objects are $0$
	and $1$ and non identity morphisms $\iota$ and $\iota^{-1}$.
	
	As similar to the homotopies of continuous functions the homotopy
	of functors are defined as follows \cite{Br1}.
	
	\begin{Def}\label{Defhomfunctor} Let $f,g\colon \mathcal{C}\rightarrow
		\mathcal{D}$ be functors. These functors are
		called {\em homotopic} and written $f\simeq g$ if there is a
		functor $F\colon \mathcal{C}\times \mathcal{J}\rightarrow
		\mathcal{D}$ such that $F(-,0)=f$ and $F(-,1)=g$ \end{Def}
	\begin{Prop}\label{Prophom}{\em  \cite{Br1} If the maps  $f,g\colon X\rightarrow
			Y$ are homotopic, then the induced
			fundamental groupoid functors
			$\pi_1f,\pi_1g\colon \pi_1X \rightarrow \pi_1Y$ are
			homotopic.}\end{Prop}
	
	\begin{Def} Let  $f,g\colon \mathcal{C}\rightarrow \mathcal{D}$  be two functors.
		We call $f$ and $g$ are {\em naturally isomorph} if there exists a natural
		isomorphism $\sigma\colon f\rightarrow g$.\end{Def}
	
	\begin{Prop} The functors  $f,g\colon \mathcal{C}\rightarrow \mathcal{D}$ are
		homotopic in the sense of Definition \ref{Defhomfunctor} if and only if they are
		naturally isomorph.\end{Prop}
	\begin{Prf} If the functors  $f,g\colon \mathcal{C}\rightarrow
		\mathcal{D}$ are homotopic  there is a functor $F\colon \mathcal{C}\times
		\mathcal{J}\rightarrow \mathcal{D}$ such that   $F(-,0)=f$ and $F(-,1)=g$. Since
		$(1_x,\iota)\colon (x,0)\rightarrow (x,1)$ is an isomorphism in $\mathcal{C}\times \mathcal{J}$ the morphism  $F(1_x,\iota)\colon F(x,0)\rightarrow F(x,1)$ is an isomorphism in $\mathcal{D}$ where $F(x,0)=f(x)$ and $F(x,1)=g(x)$.  We now define a natural  transformation  $\sigma\colon f\rightarrow g$ by
		$\sigma(x)=F(x,\iota)\colon f(x)\rightarrow g(x)$
		for  $x\in O_\mathcal{C}$.  We now prove that for a morphism $\alpha\colon x\rightarrow y$ in $\mathcal{C}$ the diagram
		\begin{equation*}
			\xymatrix{
				f(x) \ar[d]_{f(\alpha)} \ar[r]^{\sigma{(x)}} & g(x)  \ar[d]^{g(\alpha)} \\
				f(y) \ar[r]_{\sigma{(y)}} & g(y) }
		\end{equation*}
		is commutative. For this we show that the diagram
		\begin{equation*}
			\xymatrix{
				F(x,0) \ar[d]_{F(\alpha,0)} \ar[r]^{F(1_x,\iota)} & F(x,1)  \ar[d]^{F(\alpha,1)} \\
				F(y,0) \ar[r]_{F(1_y,\iota)} & F(y,1) }
		\end{equation*}
		is commutative. Since $F$ is a functor \begin{align*}
			F(\alpha,1)\circ F(1_x,\iota)& =F((\alpha,1)\circ(1_x,\iota))\\
			& =F(\alpha\circ 1_x,1\circ\iota)=F(\alpha,1)
		\end{align*}
		and
		\begin{align*}
			F(1_y,\iota)\circ F(\alpha,0) & =F((1_y,\iota)\circ(\alpha,0))\\
			& =F(1_y\circ\alpha,\iota\circ 0)=F(\alpha,1)
		\end{align*}
		and therefore the latter diagram is  commutative. Therefore the functors $f$ and $g$ are naturally isomorph.

		Conversely let   the functors  $f,g\colon \mathcal{C}\rightarrow \mathcal{D}$ be
		naturally isomorph.  So there is a natural transformation
		$\sigma\colon f\rightarrow g$  such that $\sigma_{x}\colon f(x)\rightarrow g(x)$
		is an isomorphism for each $x\in O_\mathcal{C}$ and
		for $x,y\in O_\mathcal{C}$  and  $\alpha\in \mathcal{C}(x,y)$ the following
		diagram is commutative
		\begin{equation*}
			\xymatrix{
				x \ar[d]_{\alpha} & f(x) \ar[d]_{f(\alpha)} \ar[r]^{\sigma{(x)}}    & g(x)  \ar[d]^{g(\alpha)} \\
				y                 & f(y) \ar[r]_{\sigma{(y)}}                       & g(y) }
		\end{equation*}
		We now define a homotopy of functors  $F\colon
		\mathcal{C}\times \mathcal{J}\rightarrow \mathcal{D}$ as follows:
		Define $F$ on objects by $F(x,0)=f(x)$ and
		$F(x,1)=g(x)$ for $x\in O_\mathcal{C}$. For   $x,y,z\in O_\mathcal{C}$ consider the following
		diagram of the morphisms in $\mathcal{C}\times \mathcal{J}$.
		\begin{equation*}
			\xymatrix{
				(x,0) \ar[d]_{(\alpha,\iota)} \ar[r]^{(\alpha,0)}    & (y,0)   \\
				(y,1)  & (x,1) \ar[l]^{(\alpha,1)} \ar[u]_{(\alpha,\iota^{-1})} }
		\end{equation*}
		Define  $F$ on these morphisms as follow:
		\begin{align*}
			F(\alpha,0)              &=f(\alpha)\\
			F(\alpha,1)              &=g(\alpha)\\
			F(\alpha,\iota)          &=g(\alpha)\circ\sigma_{x} \\
			F(\alpha,\iota^{-1})     &=f(\alpha)\circ(\sigma_{x})^{-1} .\end{align*}
		In this way a functor  $F\colon \mathcal{C}\times \mathcal{J}\rightarrow \mathcal{D}$
		is defined  such that $F(-,0)=f$ ve $F(-,1)=g$.  Therefore the functors $f$
		and $g$ are homotopic.
	\end{Prf}
	
	A {\em group-groupoid} which is also known  as {\em 2-group} in literature is a group object in the category of groupoids.  The formal definition of a group-groupoid is given in  \cite{BS} under the name {\em G-groupoid} as follows:
	
	\begin{Def}\label{Defgroup-groupoid} {\em  A {\em group-groupoid} $G$ is a groupoid endowed with a group structure such that the  following maps, which are called respectively product, inverse and unit   are the morphisms of groupoids:
			
			(i) $m\colon G\times G\rightarrow G$, $(g,h)\mapsto gh$;
			
			(ii) $u\colon G\rightarrow G$, $g\mapsto \overline{g}$;
			
			(iii) $e\colon \{\star\}\rightarrow G$, where $\{\star\}$ is singleton.}
	\end{Def}
	
	Here note that the group axioms can be stated as:
	
	1.  $m(1\times m)= m(m\times 1)$
	
	2.  $m~i_1= 1_G= m~i_2$
	
	3.  $m(1,u)= m(u,1)=e$
	\\ where $i_1,i_2\colon
	G\rightarrow G\times G $ are injections defined by
	$i_1(a)=(a,e)$ and $i_2(a)=(e,a)$, and  $e\colon G\rightarrow
	G$ is the constant map at $e$.

	In the definition of group-groupoid if we take these functors to be
	homotopic rather than equal, we obtain a kind of definition of
	categorical group.  There are various forms of  definitions of
	categorical group in the literature (see \cite{Car} and \cite{Mac})
	and we will use the following one with some weak conditions.
	\begin{Def}
		{\em Let $\mathcal{G}$ be a groupoid.   Let $\otimes\colon
			\mathcal{G}\times \mathcal{G}\rightarrow \mathcal{G}$ and $u\colon
			\mathcal{G}\rightarrow \mathcal{G},a\mapsto \overline{a}$ be
			functors called respectively product and inverse. Let  $e\in
			O_\mathcal{G}$ be  an object. If the following conditions are
			satisfied then $\mathcal{G}=(\mathcal{G},\otimes,u,e)$ is called a
			{\em categorical group}.
			\begin{enumerate}
				\item The functors $\otimes(1\times \otimes),\otimes(1\times \otimes)\colon \mathcal{G}\times \mathcal{G}
				\times \mathcal{G}\rightarrow \mathcal{G}$
				are homotopic.
				\item  The functors  $e \otimes 1, 1\otimes e\colon \mathcal{G}\rightarrow \mathcal{G}$
				defined by $(e\otimes 1)(a)=e \otimes a$  and $(1\times e)(a)=a \otimes e $ for $a\in \mathcal{G}$   are homotopic to the identity functor $\mathcal{G}\rightarrow \mathcal{G}$.
				\item  the functors $\otimes(1,u), \otimes(u,1)\colon \mathcal{C}\rightarrow \mathcal{C}$
				defined by $\otimes(1,u)(a)=a\otimes u(a)$ and $\otimes(u,1)(a)=u(a)\otimes a$ are homotopic to the constant functor $e\colon \mathcal{C}\rightarrow \mathcal{C}$ .
			\end{enumerate}
		}\end{Def}
		
		In this definition if  these  functors are equal rather than
		homotopic,  then the categorical group is called a {\em strict
			categorical group} which is also called group-groupoid or 2-groups.

		Note that the product  $\otimes\colon \mathcal{G}\times
		\mathcal{G}\rightarrow \mathcal{G}$ is a functor if and only if
		\[(b\circ a)\otimes (d\circ c)=(b\otimes d)\circ (a\otimes c )\]
		for $a,b,c,d\in \mathcal{G}$ whenever the compositions  $b\circ a$
		and $d\circ c$ are defined.  Since $u\colon \mathcal{G}\rightarrow
		\mathcal{G},a\mapsto \overline{a}$  is a functor when the groupoid
		composition  $b\circ a$ is defined $\overline{b\circ
			a}=\overline{b}\circ \overline{a}$  and
		$\overline{1_x}=1_{\overline{x}}$ for $x\in O_\mathcal{G}$.

		\begin{Prop}\label{fundfunctor} If $(X,x_0)$ is  an $H$-group, then the
			fundamental group $\pi_1 X$ is a categorical group. \end{Prop}
		\begin{Prf} Since $(X,x_0)$ is an $H$-group there are continuous maps
			
			\[m \colon (X\times X,(x_0,x_0))\rightarrow (X,x_0), (x,y)\mapsto x\circ y\]
			\[n\colon (X,x_0)\rightarrow (X,x_0), x\mapsto
			\overline{x}\]
			and the following  homotopies of the maps
			
			(i) $m(1_X\times m)\simeq m(m\times 1_X)$;
			
			(ii) $m~ i_1\simeq 1_X\simeq m~i_2$;
			
			(iii) $m(1_X,n)\simeq c\simeq m(n,1_X)$
			\\ where
			$i_1,i_2\colon X\rightarrow X\times X $ are injection defined by
			$i_1(x)=(x,x_0)$,  $i_2(x)=(x_0,x)$ and  $c\colon X\rightarrow X$ is the
			constant map at $x_0$.  From $m$  and $n$  we have the following
			induced functors
			\[\widetilde{m}=\pi_1m \colon \pi_1X\times \pi_1X\rightarrow \pi_1X\] and
			\[\widetilde{n}=\pi_1n \colon \pi_1X\rightarrow \pi_1X.\]
			By Proposition \ref{Prophom} from the above homotopies (i), (ii)
			and (iii)  , the following homotopies of the functors are obtained
			
			(i) $\pi_1m(1\times \pi_1m)\simeq \pi_1m(\pi_1m\times 1)$
			
			(ii) $\pi_1m~ \pi_1 i_1\simeq 1_{\pi_1 X}\simeq \pi_1m~\pi_1 i_2$
			
			(iii) $\pi_1m(1_{\pi_1 X},\pi_1n)\simeq \pi_1 c\simeq \pi_1m(\pi_1
			n,1_{\pi_1 X )}$.
			
			Therefore  $\pi_1X$ becomes  a categorical group.\end{Prf}

		\begin{Def}{\em  Let $\mathcal{G}$ and $\mathcal{H}$ be  two categorical groups.  A morphism of
				categorical groups is a morphism  $f\colon \mathcal{G}\rightarrow \mathcal{H}$ of groupoids  such that $f(a\otimes
				b)=f(a)\otimes f(b)$ for $a,b\in G$. }\end{Def}
		
		So we have  a  category denoted by $Catgroups$ of categorical groups.
		
		\begin{Prop}\label{fundfunctor1}{\em If $p\colon
				(\widetilde{X},\widetilde{x_0})\rightarrow (X,x_0)$
				is a morphism of  $H$-groups, then the induced map $\pi_{1}p\colon
				\pi_{1}(\widetilde{X})\rightarrow \pi_1(X)$ is a morphism of
				categorical groups.}\end{Prop}
		
		\begin{Exam}{\em  Let  $(X,x_0)$  be an $H$-group. Then we have  a slice category
				\[ HGpCov/(X,x_0) \] of $H$-group morphisms
				$f\colon(\widetilde{X},\widetilde{x}_0)\rightarrow (X,x_0)$ which
				are  covering maps on the underlying spaces. Hence a morphism from
				$f\colon (\widetilde{X},\widetilde{x}_0)\rightarrow (X,x_0)$ to
				$g\colon (\widetilde{Y},\widetilde{y}_0)\rightarrow (X,x_0)$ is a
				continuous map $p \colon
				(\widetilde{X},\widetilde{x}_0)\rightarrow
				(\widetilde{Y},\widetilde{y}_0)$ which becomes also  a covering
				map such that $f=gp$.
				
				Similarly we  have an other  slice category $CatGpCov/\pi_1 X$ of  categorical group morphisms  $p\colon \widetilde{G}\rightarrow \pi_1(X)$ which are
				covering morphisms on underlying groupoids. }\end{Exam}
		
		\begin{Theo} {\em  Let $(X,x_0)$ be an $H$-group such that the underlying space
				has a simply connected cover.
				Then the categories  $HGpCov/(X,x_0)$ and  $CatGpCov/\pi_1 X$ are
				equivalent. }\end{Theo}
		\begin{Prf}  Let  $p\colon (\widetilde{X},\widetilde{x}_0)\rightarrow
			(X,x_0)$ be  a morphism of $H$-groups which is a covering map on the spaces. Then by Proposition
			\ref{fundfunctor1} the induced  morphism $\pi_1 p\colon
			\pi_1\widetilde{X}\rightarrow \pi_1X$ is a morphism of categorical groups which
			is a covering morphism of underlying groupoids. So in this way we have a
			functor
			\[\pi_1\colon HGpCov/(X,x_0)\rightarrow CatGpCov/\pi_1 X.\]
			Conversely we define a functor  \[\eta\colon CatGpCov/\pi_1
			X\rightarrow HGpCov/(X,x_0)\] as follows:
			
			Let $q\colon \widetilde{G}\rightarrow
			\pi_1X$ be a morphism of categorical groups which is a covering
			morphism on the underlying groupoids. Then by 9.5.5 of \cite{Br1}
			there is a topology on $\widetilde{X}=O_{\widetilde{G}}$ and an
			isomorphism $\alpha\colon \widetilde{G}\rightarrow
			\pi_1(\widetilde{X})$ such that $p=O_q\colon
			(\widetilde{X},\widetilde{x}_0)\rightarrow (X,x_{0})$ is a
			covering map and $q=\pi_1(p)\circ\alpha$.  Hence categorical group
			structure on $\widetilde{G}$ transports via $\alpha$  to
			$\pi_1(\widetilde{X})$.  So we have the morphisms of groupoids
			\[\widetilde{m}\colon \pi_1(\widetilde{X})\times
			\pi_1(\widetilde{X})\longrightarrow \pi_1(\widetilde{X})\]
			\[\widetilde{n}\colon \pi_1(\widetilde{X})\longrightarrow \pi_1(\widetilde{X})\]
			such that $\pi_1(p)\circ \widetilde{m}=m\circ (\pi_1(p)\times
			\pi_1(p))$ and $n\pi_1(p)=\pi_1(p)\widetilde{n}$.  From these
			morphisms we obtain the  maps
			
			\[\widetilde{m}\colon \widetilde{X}\times \widetilde{X}\longrightarrow
			\widetilde{X}\]
			\[\widetilde{n}\colon \widetilde{X}\longrightarrow \widetilde{X}\]
			
			Since $(X,x_0)$ is an $H$-group with  the maps \[m \colon (X\times
			X,(x_0,x_0))\rightarrow (X,x_0)\]
			\[n\colon (X,x_0)\rightarrow (X,x_0)\]
			we have the following  homotopies of pointed maps
			
			(i)  $m(1_X\times m)\simeq m(m\times 1_X)$
			
			(ii) $m~i_1\simeq 1_X\simeq m~i_2$
			
			(iii) $m(1_X,n)\simeq c\simeq m(n,1_X)$.
			
			Then  by  Theorem \ref{homotopy covering theorem} we have the following homotopies
			
			(i)  $\widetilde{m}(1_{\widetilde{X}}\times \widetilde{m})\simeq
			\widetilde{m}(\widetilde{m}\times 1_{\widetilde{X}})$
			
			(ii) $\widetilde{m}~i_1\simeq 1_{\widetilde{X}}\simeq \widetilde{m}~i_2$
			
			(iii) $\widetilde{m}(1_{\widetilde{X}},n)\simeq c\simeq \widetilde{m}(n,1_{\widetilde{X}})$.

			Therefore $(\widetilde{X},\widetilde{x_0})$ is an $H$-group and
			$q=Ob(p)\colon (\widetilde{X},\widetilde{x}_0)\rightarrow (X,x_0)$
			is a covering morphism of $H$-groups.
			
			If $p\colon (\widetilde{X},\widetilde{x}_0) \rightarrow (X,x_0)$
			is a covering map on underlying spaces, then by  9.5.5 of
			\cite{Br1} the topology of $\widetilde{X}$ is that of $X$ lifted
			by the covering morphism $\pi_1 p\colon
			\pi_1\widetilde{X}\rightarrow \pi_1 X$ and so $\eta\pi_1=1$.
			Further if $q\colon \widetilde{G}\rightarrow \pi_1X$ is a morphism
			of categorical groups, then for the lifted topology on
			$\widetilde{X}$,   $\widetilde{G}$ is isomorph to
			$\pi_1\widetilde{X}$ and so  $\eta \pi_1\simeq 1$. Therefore these
			functors give an equivalence of the categories.
		\end{Prf}
		
		\begin{Def}{\em Let $\mathcal{G}$ be a categorical group,  $e\in O_\mathcal{G}$ the  base point
				and let  $\wtilde G$ be just a groupoid. Suppose  $p\colon\wtilde
				G\rightarrow \mathcal{G}$ is  a covering morphism of groupoids and
				$\wtilde{e}\in O_{\wtilde G}$ such that $p(\wtilde{e})=e$.
				We say that the categorical group structure of $\mathcal{G}$  {\em lifts} to $\wtilde G$
				if there exists a categorical group structure on $\wtilde G$ with the
				base point  $\wtilde{e}\in O_{\wtilde G}$ such that
				$p\colon\wtilde G\rightarrow \mathcal{G}$  is a morphism of  categorical groups }\end{Def}
		
		In the following proposition we prove that the liftings of
		homotopic functors are also homotopic.
		\begin{Prop} \label{TheoLiftedHgp1}
			Let $p\colon (\widetilde{G},\widetilde{x}) \rightarrow (G,x)$ be a
			covering morphism of groupoids. Suppose that $K$ is a 1-connected
			groupoid, i.e, for each $x,y\in O_K$, $K(x,y)$ has only one
			morphism. Let $f,g\colon (K,z)\rightarrow (G,x)$ be the morphisms of
			groupoids such that $f$ and $g$ are homotopics. Let
			$\widetilde{f}$  and $\widetilde{g}$ be the liftings of $f$ and
			$g$ respectively. Then $\widetilde{f}$ and $\widetilde{g}$ are
			also homotopic.
		\end{Prop}
		\begin{Prf} Since the functors $f$ and $g$ are homotopic, there is
			a functor $F\colon K\times \mathcal{J}\rightarrow G$ such that
			$F(-,0)=f$ and $F(-,1)=g$.  Since $K$ is a 1-connected groupoid by
			Theorem \ref{Theolifted} there is a functor $\widetilde{F}\colon (K
			\times \mathcal{J}, (z,0))\rightarrow (\widetilde{G},\widetilde{x})$
			such that $p\widetilde{F}=F$.  Hence
			$p\widetilde{F}(-,0)=F(-,0)=f$ and $p\widetilde{F}(-,1)=F(-,1)=g$.
			So by the uniqueness of the liftings we have that
			$\widetilde{F}(-,0)=\widetilde{f}$ and
			$\widetilde{F}(-,1)=\widetilde{g}$. Therefore $\widetilde{f}$ and
			$\widetilde{g}$ are homotopic.
		\end{Prf}
		
		\begin{Theo} \label{TheoLiftedHgp}
			Let $\wtilde G$ be a 1-connected groupoid and  $\mathcal{G}$  a categorical group. Suppose that  $p\colon\wtilde G\rightarrow \mathcal{G}$ is a covering morphism  on the underlying groupoids. Let $e\in O_G$ be the
			base point of $\mathcal{G}$ and $\wtilde {e}\in O_{\wtilde G}$ such that
			$p(\wtilde e)=e$. Then the categorical group structure of $\mathcal{G}$ lifts to $\wtilde G$.
		\end{Theo}
		\begin{Prf} Since $\mathcal{G}$ is categorical group, it has the following functors
			
			$\otimes\colon \mathcal{G}\times \mathcal{G}\rightarrow \mathcal{G},(a,b)\mapsto a\otimes b$
			
			$ u\colon G\rightarrow G,  a\mapsto \overline{a}$
			and
			
			$e\colon  \{\star\}\rightarrow G$\\
			such that the following reduced functors are homotopic:
			
			(i) $\otimes(1\times \otimes)\simeq \otimes(\otimes\times 1)$;
			
			(ii) $\otimes i_1\simeq  \otimes i_2\simeq 1_G$;
			
			(iii) $ \otimes(1,n)\simeq c\simeq \otimes(n,1)$.\\
			Since $\widetilde{G}$ is a 1- connected groupoid by  Theorem
			\ref{Theolifted} the functors $\otimes$ and $u$ lift respectively to the
			morphisms of groupoids
			\[\widetilde{\otimes}\colon (\widetilde{G}\times
			\widetilde{G},(\widetilde{e},\widetilde{e}))\rightarrow
			(\widetilde{G},\widetilde{e})\]
			
			and  \[\widetilde{u}\colon  (\widetilde{G},\widetilde{e})\rightarrow
			(\widetilde{G}, \widetilde{e})\]
			By Proposition  \ref{TheoLiftedHgp1} we have the following
			homotopies of the functors:
			
			(i) $\widetilde{\otimes}(1\times \widetilde{\otimes})\simeq
			\widetilde{\otimes}(\widetilde{\otimes}\times 1)$;
			
			(ii) $\widetilde{\otimes}~ i_1\simeq  \widetilde{\otimes}~i_2\simeq
			1_{\widetilde{G}}$;
			
			(iii) $ \widetilde{\otimes}(1,\widetilde{u})\simeq c\simeq
			\widetilde{\otimes}(\widetilde{u},1)$.
		\end{Prf}
		
		As a result of Theorem \ref{TheoLiftedHgp} we obtain the following
		corollary.
		\begin{Cor}\label{Corlifting}{\em Let $(X,x_0)$ be an $H$-group and  $p\colon
				(\widetilde{X},\widetilde{x}_0)\rightarrow (X,x_0)$  a covering
				map. If $\widetilde{X}$  is a simply connected topological space,
				then $H$-group structure of $(X,x_0)$ lifts to
				$(\widetilde{X},\widetilde{x}_0)$, i.e,
				$(\widetilde{X},\widetilde{x}_0)$     is an $H$-group and $p\colon
				(\widetilde{X},\widetilde{x}_0)\rightarrow (X,x_0)$ is a morphism
				of $H$-groups.}\end{Cor}
		\begin{Prf}
			Since $p\colon (\widetilde{X},\widetilde{x}_0)\rightarrow (X,x_0)$
			is a  covering map,  the induced morphism $\pi_1 p\colon\pi_1
			\wtilde X\rightarrow \pi_1X$ is a covering morphism of groupoids.
			Since $(X,x_0)$ is an $H$-group by Proposition \ref{fundfunctor}
			$\pi_1X$ is a categorical group and since $\widetilde{X}$ is
			simpy connected the fundamental groupoid $\pi_1\widetilde{X}$ is a
			simply connected groupoid. So by Proposition \ref{TheoLiftedHgp1}
			the categorical group structure of $\pi_1X$ lifts to $\pi_1 \wtilde X$. So
			we have the groupoid morphisms
			\[\widetilde{m}\colon \pi_1(\widetilde{X})\times
			\pi_1(\widetilde{X})\longrightarrow \pi_1(\widetilde{X})\]
			\[\widetilde{n}\colon \pi_1(\widetilde{X})\longrightarrow \pi_1(\widetilde{X})\]
			such that $\pi_1(p)\circ \widetilde{m}=m\circ (\pi_1(p)\times
			\pi_1(p))$ and $n\pi_1(p)=\pi_1(p)\widetilde{n}$ and so the maps
			\[\widetilde{m}\colon \widetilde{X}\times \widetilde{X}\longrightarrow
			\widetilde{X}\]
			\[\widetilde{n}\colon \widetilde{X}\longrightarrow \widetilde{X}\]
			
			Since $(X,x_0)$ is an $H$-group, we have the  homotopies
			
			(i)  $m(1_X\times m)\simeq m(m\times 1_X)$
			
			(ii) $m~i_1\simeq 1_X\simeq m~i_2$
			
			(iii) $m(1_X,n)\simeq c\simeq m(n,1_X)$
			
			and so  by  Theorem \ref{homotopy covering theorem}  the
			homotopies
			
			(i)  $\widetilde{m}(1_X\times \widetilde{m})\simeq
			\widetilde{m}(\widetilde{m}\times 1_X)$
			
			(ii) $\widetilde{m}~i_1\simeq 1_X\simeq \widetilde{m}~i_2$
			
			(iii) $\widetilde{m}(1_X,n)\simeq c\simeq \widetilde{m}(n,1_X)$.

			Therefore $(\widetilde{X},\widetilde{x_0})$ is an $H$-group and
			$q=Ob(p)\colon (\widetilde{X},\widetilde{x}_0)\rightarrow (X,x_0)$
			is a covering morphism of $H$-groups.

			Therefore $(\widetilde{X}, {\widetilde{x}_0})$ becomes an
			$H$-group as required.
		\end{Prf}

		{\bf Acknowledgement}: We would like to thank R. Brown and
		T.Datuashvili for their useful comments.

	\end{document}